\documentclass{article}
\usepackage{amssymb}
\usepackage{amsmath}

\newtheorem{theorem}{Theorem}

\newtheorem{corollary}{Corollary}

\newtheorem{lemma}{Lemma}

\newtheorem{remark}{Remark}

\begin{document}

\title{\textbf{On estimates of deviation of conjugate functions\ from matrix
operators of their Fourier series by some expressions with }$r$ \textbf{-
differences of the entries}}
\author{\textbf{W\l odzimierz \L enski \ and Bogdan Szal} \\
University of Zielona G\'{o}ra\\
Faculty of Mathematics, Computer Science and Econometrics\\
65-516 Zielona G\'{o}ra, ul. Szafrana 4a, Poland\\
W.Lenski@wmie.uz.zgora.pl , B.Szal @wmie.uz.zgora.pl}
\date{}
\maketitle

\begin{abstract}
We extend the results of the authors from [Abstract and Applied Analysis,
Volume 2016, Article ID 9712878] to the case conjugate Fourier series .

\ \ \ \ \ \ \ \ \ \ \ \ \ \ \ \ \ \ \ \ 


\ \ \ \ \ \ \ \ \ \ \ \ \ \ \ \ \ \ \ 

\end{abstract}

\section{Introduction}

Let $X=L^{p}$ or $X=C$, where $L^{p}\ (1\leq p\leq \infty )\;\left[ C\right] 
$ be the class of all $2\pi $--periodic real--valued functions, integrable
in the Lebesgue sense with $p$--th power when $p\geq 1$ and essentially
bounded when $p=\infty $ $\left[ \text{continuous}\right] $ over $Q=$ $[-\pi
,\pi ]$ with the norm%
\begin{eqnarray*}
\Vert f\Vert _{_{L^{p}}} &:&=\Vert f(\cdot )\Vert _{_{L^{p}}}=\left\{ 
\begin{array}{c}
\left( \int_{_{_{Q}}}\mid f(t)\mid ^{p}dt\right) ^{1/p}\text{ \ \ when \ }%
1\leq p<\infty , \\ 
ess\sup\limits_{t\in Q}\mid f(t)\mid \text{ \ when \ }p=\infty ,%
\end{array}%
\right. \\
\Vert f\Vert _{_{C}} &:&=\Vert f(\cdot )\Vert _{_{C}}=\sup\limits_{t\in
Q}\mid f(t)\mid
\end{eqnarray*}%
and consider the trigonometric Fourier series 
\begin{equation*}
Sf(x):=\frac{a_{0}(f)}{2}+\sum_{\nu =1}^{\infty }(a_{\nu }(f)\cos \nu
x+b_{\nu }(f)\sin \nu x)
\end{equation*}%
with the partial sums\ $S_{k}f$ and the conjugate one 
\begin{equation*}
\widetilde{S}f(x):=\sum_{\nu =1}^{\infty }(a_{\nu }(f)\sin \nu x-b_{\nu
}(f)\cos \nu x)
\end{equation*}%
with the partial sums $\widetilde{S}_{k}f$. We know that if $f\in L,$ then 
\begin{equation*}
\widetilde{f}\left( x\right) :=-\frac{1}{\pi }\int_{0}^{\pi }\psi _{x}\left(
t\right) \frac{1}{2}\cot \frac{t}{2}dt=\lim_{\epsilon \rightarrow 0^{+}}%
\widetilde{f}\left( x,\epsilon \right) ,
\end{equation*}%
where 
\begin{equation*}
\widetilde{f}\left( x,\epsilon \right) :=-\frac{1}{\pi }\int_{\epsilon
}^{\pi }\psi _{x}\left( t\right) \frac{1}{2}\cot \frac{t}{2}dt,
\end{equation*}%
with $\psi _{x}\left( t\right) :=f\left( x+t\right) -f\left( x-t\right) ,$
exists for almost all \ $x$ \cite[Th.(3.1)IV]{Z}.

Let $A:=\left( a_{n,k}\right) $ be an infinite matrix of real numbers such
that%
\begin{equation*}
a_{n,k}\geq 0\text{ when \ }k,n=0,1,2,...\text{, \ }\lim_{n\rightarrow
\infty }a_{n,k}=0\text{\ and }\sum_{k=0}^{\infty }a_{n,k}=1.
\end{equation*}%
We will use the notation $A_{n,r}=\sum_{k=0}^{\infty }\left\vert
a_{n,k}-a_{n,k+r}\right\vert ,$ for $r\in 
\mathbb{N}
.$

The\ $A-$transformation$\ $of \ $S_{k}f$ and of $\widetilde{S}f$\ be given,
by a matrix convention, as follows%
\begin{equation*}
\left( 
\begin{array}{c}
T_{n,A}^{\text{ }}f\left( x\right) \\ 
\widetilde{T}_{n,A}^{\text{ }}f\left( x\right)%
\end{array}%
\right) :=\sum_{k=0}^{\infty }a_{n,k}\left( 
\begin{array}{c}
S_{k}f\left( x\right) \\ 
\widetilde{S}_{k}f\left( x\right)%
\end{array}%
\right) \text{ \ \ \ }\left( n=0,1,2,...\right) .
\end{equation*}%
In this paper, we study the upper bounds of $\left\Vert \widetilde{T}_{n,A}^{%
\text{ }}f-\widetilde{f}\right\Vert _{X}$ and $\left\Vert \widetilde{T}%
_{n,A}^{\text{ }}f\left( \cdot \right) -\widetilde{f}\left( \cdot ,\epsilon
\right) \right\Vert _{X}$\ by the second conjugate modulus of continuity of $%
f$ in the space $X$ defined by the formula%
\begin{equation*}
\widetilde{\omega }_{2}\left( f,\delta \right) _{X}=\sup_{0<t\leq \delta
}\Vert \psi _{\cdot }\left( t\right) \Vert _{X}.
\end{equation*}%
We will also used the second classical modulus of continuity of $f$ in the
space $X$ defined by the formula $\omega _{2}\left( f,\delta \right)
_{X}:=\sup_{0\leq \left\vert t\right\vert \leq \delta }\left\Vert \varphi
_{.}\left( t\right) \right\Vert _{X},$ where $\varphi _{x}\left( t\right)
:=f\left( x+t\right) +f\left( x-t\right) -2f\left( x\right) .$

We will consider a function $\omega $ of modulus of continuity type on the
interval $[0,2\pi ],$ i.e. a nondecreasing continuous function having the
following properties:\ $\omega \left( 0\right) =0,$\ $\omega \left( \delta
_{1}+\delta _{2}\right) \leq \omega \left( \delta _{1}\right) +\omega \left(
\delta _{2}\right) $\ for any\ $0\leq \delta _{1}\leq \delta _{2}\leq \delta
_{1}+\delta _{2}\leq 2\pi $.

The deviation $T_{n,A}^{\text{ }}f-f$\ was estimated in the paper \cite%
{LWBSZ} ( see also \cite[Theorems 3.4, p. 290]{XK} and \cite{WY}) as follows:

\textbf{Theorem A.} \textit{Let }$f\in X_{\omega }=\left\{ f\in X:\omega
_{2}\left( f,\delta \right) _{X}=O\left( \omega \left( \delta \right)
\right) \text{ \textit{when\ }}\delta \in \left[ 0,2\pi \right] \right\} $%
\textit{\ and }$r\in 
\mathbb{N}
.$\textit{\ Then,}%
\begin{equation*}
\left\Vert T_{n,A}^{\text{ }}f-f\right\Vert _{X}=O\left( H\left( \frac{\pi }{%
n+1}\right) \left( \frac{\pi }{n+1}+A_{n,r}\right) \right) ,
\end{equation*}%
\textit{where a function of modulus of continuity type }$\omega $ \textit{%
satisfies the condition\ }%
\begin{equation}
\int_{u}^{\pi }t^{-2}\omega \left( t\right) dt=O\left( H\left( u\right)
\right) \text{ \textit{when}\ }u\in \left[ 0,\pi \right] ,  \label{111}
\end{equation}%
\textit{with }$H(u)\geq 0$\textit{, such that}%
\begin{equation}
\int_{0}^{u}H\left( t\right) dt=O\left( uH\left( u\right) \right) \text{ 
\textit{when}\ }u\in \left[ 0,\pi \right] .  \label{112}
\end{equation}%
\textit{Additionally, if }%
\begin{equation}
\left[ \sum_{l=0}^{n}\sum_{k=l}^{r+l-1}a_{n,k}\right] ^{-1}=O\left( 1\right)
,  \label{113}
\end{equation}%
\textit{then}%
\begin{equation*}
\left\Vert T_{n,A}^{\text{ }}f-f\right\Vert _{X}=O\left( H\left( \frac{\pi }{%
n+1}\right) A_{n,r}\right)
\end{equation*}%
\textit{but if }%
\begin{equation}
\sum_{k=0}^{\infty }\left( k+1\right) a_{n,k}=O\left( n+1\right) ,
\label{114}
\end{equation}%
\textit{then}%
\begin{equation*}
\left\Vert T_{n,A}^{\text{ }}f-f\right\Vert _{X}=O\left( \omega \left( \frac{%
\pi }{n+1}\right) +H\left( \frac{\pi }{n+1}\right) A_{n,r}\right) .
\end{equation*}%
\textbf{Theorem B.} \textit{If }$f\in X$\textit{\ and a matrix }$A$\textit{\
is such that }$\left( \ref{114}\right) $\textit{\ holds, then for }$r\in 
\mathbb{N}
$\textit{\ }%
\begin{eqnarray*}
\left\Vert T_{n,A}^{\text{ }}f-f\right\Vert _{X} &=&O\left( \omega
_{2}\left( f,\frac{\pi }{n+1}\right) _{X}+\sum_{\mu =1}^{n}\frac{\omega
_{2}\left( f,\frac{\pi }{\mu }\right) _{X}}{\mu }\sum_{k=0}^{\mu
+1}a_{n,k}\right. \\
&&+\left. \sum_{\mu =1}^{n}\omega _{2}\left( f,\frac{\pi }{\mu }\right)
_{X}\sum_{k=\mu }^{\infty }\left\vert a_{n,k}-a_{n,k+r}\right\vert \right) .
\end{eqnarray*}%
\textit{\ }

From our theorems we also derived a corollary for a matrix $A$ satisfying
the following condition $\sum_{k=m}^{\infty }\left\vert
a_{n,k}-a_{n,k+r}\right\vert \leq K\sum_{k=n/c}^{\infty }\frac{a_{n,k}}{k}.$

\section{Statement of the results}

Let $X_{\omega }=\left\{ f\in X:\widetilde{\omega }_{2}\left( f,\delta
\right) _{X}=O\left( \omega \left( \delta \right) \right) \text{ when\ }%
\delta \in \left[ 0,2\pi \right] \right\} .$ We present the estimates of the
quantities $\left\Vert \widetilde{T}_{n,A}^{\text{ }}f\left( \cdot \right) -%
\widetilde{f}\left( \cdot \right) \right\Vert _{X}$ and $\left\Vert 
\widetilde{T}_{n,A}^{\text{ }}f\left( \cdot \right) -\widetilde{f}\left(
\cdot ,\epsilon \right) \right\Vert _{X}$ simultaneously. Finally, we give a
corollary and a remark.

\begin{theorem}
If $f\in X_{\omega },$ where $\omega $ satisfies the condition $\left( \ref%
{111}\right) $ such that $\left( \ref{112}\right) $ holds and $r\in 
\mathbb{N}
$, then%
\begin{equation*}
\left. 
\begin{array}{c}
\left\Vert \widetilde{T}_{n,A}^{\text{ }}f\left( \cdot \right) -\widetilde{f}%
\left( \cdot \right) \right\Vert _{X} \\ 
\left\Vert \widetilde{T}_{n,A}^{\text{ }}f\left( \cdot \right) -\widetilde{f}%
\left( \cdot ,\frac{\pi }{r\left( n+1\right) }\right) \right\Vert _{X}%
\end{array}%
\right\} =O\left( H\left( \frac{\pi }{n+1}\right) \left( \frac{\pi }{n+1}%
+A_{n,r}\right) \right) .
\end{equation*}%
Additionally, if a matrix $A$ is such that $\left( \ref{113}\right) $ is
true, then%
\begin{equation*}
\left. 
\begin{array}{c}
\left\Vert \widetilde{T}_{n,A}^{\text{ }}f\left( \cdot \right) -\widetilde{f}%
\left( \cdot \right) \right\Vert _{X} \\ 
\left\Vert \widetilde{T}_{n,A}^{\text{ }}f\left( \cdot \right) -\widetilde{f}%
\left( \cdot ,\frac{\pi }{r\left( n+1\right) }\right) \right\Vert _{X}%
\end{array}%
\right\} =O\left( H\left( \frac{\pi }{n+1}\right) A_{n,r}\right) .
\end{equation*}
\end{theorem}

\begin{theorem}
If $f\in X_{\omega },$ where $\omega $ satisfies the condition $\left( \ref%
{111}\right) $ such that $\left( \ref{112}\right) $ holds, $r\in 
\mathbb{N}
$ and a matrix $A$ is such that $\left( \ref{113}\right) $ is true, then 
\begin{equation*}
\left. 
\begin{array}{c}
\left\Vert \widetilde{T}_{n,A}^{\text{ }}f\left( \cdot \right) -\widetilde{f}%
\left( \cdot \right) \right\Vert _{X} \\ 
\left\Vert \widetilde{T}_{n,A}^{\text{ }}f\left( \cdot \right) -\widetilde{f}%
\left( \cdot ,\frac{1}{r}A_{n,r}\right) \right\Vert _{X}%
\end{array}%
\right\} =O\left( H\left( A_{n,r}\right) A_{n,r}\right) .
\end{equation*}
\end{theorem}

\begin{theorem}
If $f\in X_{\omega },$ where $\omega $ satisfies the condition $\left( \ref%
{111}\right) $ such that $\left( \ref{112}\right) $ holds and $r\in 
\mathbb{N}
$, then%
\begin{equation*}
\left. 
\begin{array}{c}
\left\Vert \widetilde{T}_{n,A}^{\text{ }}f\left( \cdot \right) -\widetilde{f}%
\left( \cdot \right) \right\Vert _{X} \\ 
\left\Vert \widetilde{T}_{n,A}^{\text{ }}f\left( \cdot \right) -\widetilde{f}%
\left( \cdot ,\frac{\pi }{r\left( n+1\right) }\right) \right\Vert _{X}%
\end{array}%
\right\} =O\left( \omega \left( \frac{\pi }{n+1}\right) +H\left( \frac{\pi }{%
n+1}\right) A_{n,r}\right) ,
\end{equation*}%
were in the case of the first estimate $\omega $ satisfies the extra
condition%
\begin{equation}
\int_{0}^{u}t^{-1}\omega \left( t\right) dt=O\left( \omega \left( u\right)
\right) \text{ when\ }u\in \left[ 0,2\pi \right] ,  \label{202}
\end{equation}%
but in the case of the second estimate a matrix $A$ is such that $\ \left( %
\ref{114}\right) $ is true.
\end{theorem}

\begin{theorem}
If $f\in X$ and $r\in 
\mathbb{N}
$, then%
\begin{eqnarray*}
\left. 
\begin{array}{c}
\left\Vert \widetilde{T}_{n,A}^{\text{ }}f\left( \cdot \right) -\widetilde{f}%
\left( \cdot \right) \right\Vert _{X} \\ 
\left\Vert \widetilde{T}_{n,A}^{\text{ }}f\left( \cdot \right) -\widetilde{f}%
\left( \cdot ,\frac{\pi }{r\left( n+1\right) }\right) \right\Vert _{X}%
\end{array}%
\right\} &=&O\left( \widetilde{\omega }_{2}\left( f,\frac{\pi }{n+1}\right)
_{X}+\sum_{\mu =1}^{n}\frac{\widetilde{\omega }_{2}\left( f,\frac{\pi }{\mu }%
\right) _{X}}{\mu }\sum_{k=0}^{\mu +1}a_{n,k}\right. \\
&&+\left. \sum_{\mu =1}^{n}\widetilde{\omega }_{2}\left( f,\frac{\pi }{\mu }%
\right) _{X}\sum_{k=\mu }^{\infty }\left\vert a_{n,k}-a_{n,k+r}\right\vert
\right) ,
\end{eqnarray*}%
were in the case of the first estimate $\widetilde{\omega }_{2}$ instead of $%
\omega $ satisfies the extra condition $\left( \ref{202}\right) $, but in
the case of the second estimate a matrix $A$ is such that $\left( \ref{114}%
\right) $ is true.
\end{theorem}

\begin{corollary}
From Theorem 4 it follows that if $f\in X_{\omega },$ where $\omega $
satisfies the condition $\left( \ref{111}\right) $ such that $\left( \ref%
{112}\right) $ is true and%
\begin{equation}
\sum_{k=m}^{\infty }\left\vert a_{n,k}-a_{n,k+r}\right\vert \leq
K\sum_{k=m/c}^{\infty }\frac{a_{n,k}}{k},  \label{200}
\end{equation}%
with some $c>1$ and $r\in 
\mathbb{N}
,$ holds, then%
\begin{equation*}
\left. 
\begin{array}{c}
\left\Vert \widetilde{T}_{n,A}^{\text{ }}f\left( \cdot \right) -\widetilde{f}%
\left( \cdot \right) \right\Vert _{X} \\ 
\left\Vert \widetilde{T}_{n,A}^{\text{ }}f\left( \cdot \right) -\widetilde{f}%
\left( \cdot ,\frac{\pi }{r\left( n+1\right) }\right) \right\Vert _{X}%
\end{array}%
\right\} =O\left( \frac{H\left( \frac{\pi }{n+1}\right) }{n+1}%
+\sum_{k=1}^{n}a_{n,k}\frac{H\left( \frac{\pi }{k+1}\right) }{k+1}\right) ,
\end{equation*}%
were in the case of the first estimate $\widetilde{\omega }_{2}$ instead of $%
\omega $ satisfies the extra condition $\left( \ref{202}\right) $, but in
the case of the second estimate a matrix $A$ is such that $\left( \ref{114}%
\right) $ is true.
\end{corollary}

\begin{remark}
We note that our extra conditions $\left( \ref{113}\right) $ and $\left( \ref%
{114}\right) $ for a lower triangular infinite matrix\ $A$ always hold.
\end{remark}

\section{Auxiliary results}

We begin\ this section by some notations from \cite{BSZal2} and \cite[%
Section 5 of Chapter II]{Z}. Let for $r=1,2,...$%
\begin{equation*}
D_{k,r}^{\circ }\left( t\right) =\frac{\sin \frac{\left( 2k+r\right) t}{2}}{%
2\sin \frac{rt}{2}}\text{, }\widetilde{D^{\circ }}_{k,r}\left( t\right) =%
\frac{\cos \frac{\left( 2k+r\right) t}{2}}{2\sin \frac{rt}{2}}\text{ and }%
\widetilde{D}_{k,r}\left( t\right) =\frac{\cos \frac{t}{2}-\cos \frac{\left(
2k+r\right) t}{2}}{2\sin \frac{rt}{2}}.
\end{equation*}

It is clear by \cite{Z} that $\widetilde{S}_{k}f\left( x\right) =-\frac{1}{%
\pi }\int_{-\pi }^{\pi }f\left( x+t\right) \widetilde{D}_{k,1}\left(
t\right) dt,$ whence%
\begin{equation*}
\widetilde{T}_{n,A}^{\text{ }}f\left( x\right) -\widetilde{f}\left( x\right)
=\frac{1}{\pi }\int_{0}^{\pi }\psi _{x}\left( t\right) \sum_{k=0}^{\infty
}a_{n,k}\widetilde{D^{\circ }}_{k,1}\left( t\right) dt
\end{equation*}%
and 
\begin{eqnarray*}
\widetilde{T}_{n,A}^{\text{ }}f\left( x\right) -\widetilde{f}\left( x,\frac{%
\pi }{r\left( n+1\right) }\right) &=&-\frac{1}{\pi }\int_{0}^{\frac{\pi }{%
r\left( n+1\right) }}\psi _{x}\left( t\right) \sum_{k=0}^{\infty }a_{n,k}%
\widetilde{D}_{k,1}\left( t\right) dt \\
&&+\frac{1}{\pi }\int_{\frac{\pi }{r\left( n+1\right) }}^{\pi }\psi
_{x}\left( t\right) \sum_{k=0}^{\infty }a_{n,k}\widetilde{D^{\circ }}%
_{k,1}\left( t\right) dt.
\end{eqnarray*}

Next, we present the known estimates and relations.

\begin{lemma}
$\cite{Z}$ If \ $0<\left\vert t\right\vert \leq \pi $ then 
\begin{equation*}
\left\vert \widetilde{D^{\circ }}_{k,1}\left( t\right) \right\vert \leq 
\frac{\pi }{2\left\vert t\right\vert },\text{ }\left\vert \widetilde{D}%
_{k,1}\left( t\right) \right\vert \leq \frac{\pi }{\left\vert t\right\vert }%
\text{\ }
\end{equation*}%
and, for any real $t,$ we have%
\begin{equation*}
\left\vert \widetilde{D}_{k,1}\left( t\right) \right\vert \leq \frac{1}{2}%
k\left( k+1\right) \left\vert t\right\vert \text{ \ and }\left\vert 
\widetilde{D}_{k,1}\left( t\right) \right\vert \leq k+1.
\end{equation*}
\end{lemma}

\begin{lemma}
$\cite{BSZal2}$\textit{\ Let }$r\in N,$\textit{\ }$l\in Z$\textit{\ and }$%
a:=(a_{n})\subset 
\mathbb{C}
$\textit{. If }$x$\textit{\ }$\neq $\textit{\ }$\frac{2l\pi }{r}$\textit{\ ,
then for every }$m\geq n$%
\begin{eqnarray*}
\sum_{k=n}^{m}a_{k}\sin kx &=&-\sum_{k=n}^{m}\left( a_{k}-a_{k+r}\right) 
\widetilde{D^{\circ }}_{k,r}\left( t\right) +\sum_{k=m+1}^{m+r}a_{k}%
\widetilde{D^{\circ }}_{k,-r}\left( t\right) -\sum_{k=n}^{n+r-1}a_{k}%
\widetilde{D^{\circ }}_{k,-r}\left( t\right) , \\
\sum_{k=n}^{m}a_{k}\cos kx &=&\sum_{k=n}^{m}\left( a_{k}-a_{k+r}\right)
D_{k,r}^{\circ }\left( t\right) -\sum_{k=m+1}^{m+r}a_{k}D_{k,-r}^{\circ
}\left( t\right) +\sum_{k=n}^{n+r-1}a_{k}D_{k,-r}^{\circ }\left( t\right) .
\end{eqnarray*}
\end{lemma}

We additionally need two estimates with a function\ of modulus of continuity
type $\omega $.

\begin{lemma}
$\cite{LWBSZ}$ If $\left( \ref{111}\right) $ and $\left( \ref{112}\right) $
hold, then, for $c\geq 1$ and $\beta >\alpha >0,$%
\begin{equation*}
\int_{\alpha }^{\beta }t^{-1}\omega \left( t\right) dt=O\left( \left( \beta
-\alpha \right) H\left( c\left( \beta -\alpha \right) \right) \right) \text{
when\ }\left( \beta -\alpha \right) \in \left[ 0,2\pi \right] .
\end{equation*}
\end{lemma}

\begin{lemma}
$\cite{LWBSZ}$ If $\left( \ref{111}\right) $ and $\left( \ref{112}\right) $
hold, then, for $b\geq 1,$%
\begin{equation*}
\int_{u}^{\pi }t^{-2}\omega \left( t\right) dt=O\left( H\left( bu\right)
\right) \text{ when\ }u\in \left[ 0,\pi \right] .
\end{equation*}
\end{lemma}

Finally, we present very useful property of such function $\omega $.

\begin{lemma}
$\cite{Z}$ A function $\omega $ of modulus of continuity type on the
interval $[0,2\pi ]$ satisfies the following condition $\delta
_{2}^{-1}\omega \left( \delta _{2}\right) \leq 2\delta _{1}^{-1}\omega
\left( \delta _{1}\right) $ for$\ \delta _{2}\geq \delta _{1}>0.$ \ 
\end{lemma}

\section{Proofs of the results}

\subsection{Proof of Theorem 1}

It is clear that for an odd $r$%
\begin{eqnarray*}
&&\left( 
\begin{array}{c}
\widetilde{T}_{n,A}^{\text{ }}f\left( x\right) -\widetilde{f}\left( x\right)
\\ 
\widetilde{T}_{n,A}^{\text{ }}f\left( x\right) -\widetilde{f}\left( x,\frac{%
\pi }{r\left( n+1\right) }\right)%
\end{array}%
\right) \\
&=&\left( 
\begin{array}{c}
+ \\ 
-%
\end{array}%
\right) \frac{1}{\pi }\int\limits_{0}^{\frac{\pi }{r\left( n+1\right) }}\psi
_{x}\left( t\right) \sum_{k=0}^{\infty }a_{n,k}\left( 
\begin{array}{c}
\widetilde{D^{\circ }}_{k,1}\left( t\right) \\ 
\widetilde{D}_{k,1}\left( t\right)%
\end{array}%
\right) dt+\frac{1}{\pi }\int\limits_{\frac{\pi }{r\left( n+1\right) }}^{%
\frac{\pi }{r}}\psi _{x}\left( t\right) \sum_{k=0}^{\infty }a_{n,k}%
\widetilde{D^{\circ }}_{k,1}\left( t\right) dt \\
&&+\frac{1}{\pi }\sum_{m=1}^{\left[ r/2\right] }\int\limits_{\frac{2m\pi }{r}%
}^{\frac{2m\pi }{r}+\frac{\pi }{r}}\psi _{x}\left( t\right)
\sum_{k=0}^{\infty }a_{n,k}\widetilde{D^{\circ }}_{k,1}\left( t\right) dt+%
\frac{1}{\pi }\sum_{m=0}^{\left[ r/2\right] -1}\int\limits_{\frac{2m\pi }{r}+%
\frac{\pi }{r}}^{\frac{2\left( m+1\right) \pi }{r}}\psi _{x}\left( t\right)
\sum_{k=0}^{\infty }a_{n,k}\widetilde{D^{\circ }}_{k,1}\left( t\right) dt \\
&=&\left( 
\begin{array}{c}
J_{1}\left( x\right) \\ 
J_{1}^{\prime }\left( x\right)%
\end{array}%
\right) +J_{2}\left( x\right) +I_{1}^{\prime \prime }\left( x\right)
+I_{2}\left( x\right)
\end{eqnarray*}%
and for an even $r$%
\begin{eqnarray*}
&&\left( 
\begin{array}{c}
\widetilde{T}_{n,A}^{\text{ }}f\left( x\right) -\widetilde{f}\left( x\right)
\\ 
\widetilde{T}_{n,A}^{\text{ }}f\left( x\right) -\widetilde{f}\left( x,\frac{%
\pi }{r\left( n+1\right) }\right)%
\end{array}%
\right) \\
&=&\left( 
\begin{array}{c}
+ \\ 
-%
\end{array}%
\right) \frac{1}{\pi }\int\limits_{0}^{\frac{\pi }{r\left( n+1\right) }}\psi
_{x}\left( t\right) \sum_{k=0}^{\infty }a_{n,k}\left( 
\begin{array}{c}
\widetilde{D^{\circ }}_{k,1}\left( t\right) \\ 
\widetilde{D}_{k,1}\left( t\right)%
\end{array}%
\right) dt+\frac{1}{\pi }\int\limits_{\frac{\pi }{r\left( n+1\right) }}^{%
\frac{\pi }{r}}\psi _{x}\left( t\right) \sum_{k=0}^{\infty }a_{n,k}%
\widetilde{D^{\circ }}_{k,1}\left( t\right) dt \\
&&+\frac{1}{\pi }\sum_{m=1}^{\left[ r/2\right] -1}\int\limits_{\frac{2m\pi }{%
r}}^{\frac{2m\pi }{r}+\frac{\pi }{r}}\psi _{x}\left( t\right)
\sum_{k=0}^{\infty }a_{n,k}\widetilde{D^{\circ }}_{k,1}\left( t\right) dt+%
\frac{1}{\pi }\sum_{m=0}^{\left[ r/2\right] -1}\int\limits_{\frac{2m\pi }{r}+%
\frac{\pi }{r}}^{\frac{2\left( m+1\right) \pi }{r}}\psi _{x}\left( t\right)
\sum_{k=0}^{\infty }a_{n,k}\widetilde{D^{\circ }}_{k,1}\left( t\right) dt \\
&=&\left( 
\begin{array}{c}
J_{1}\left( x\right) \\ 
J_{1}^{\prime }\left( x\right)%
\end{array}%
\right) +J_{2}\left( x\right) +I_{1}^{\prime \prime \prime }\left( x\right)
+I_{2}\left( x\right) .
\end{eqnarray*}%
Then,%
\begin{equation*}
\left( 
\begin{array}{c}
\left\Vert \widetilde{T}_{n,A}^{\text{ }}f\left( \cdot \right) -\widetilde{f}%
\left( \cdot \right) \right\Vert _{X} \\ 
\left\Vert \widetilde{T}_{n,A}^{\text{ }}f\left( \cdot \right) -\widetilde{f}%
\left( \cdot ,\frac{\pi }{r\left( n+1\right) }\right) \right\Vert _{X}%
\end{array}%
\right) \leq \left( 
\begin{array}{c}
\left\Vert J_{1}+J_{2}+I_{1}^{\prime \prime }\right\Vert _{X}+\left\Vert
J_{1}+J_{2}+I_{1}^{\prime \prime \prime }\right\Vert _{X} \\ 
\left\Vert J_{1}^{\prime }\right\Vert _{X}+\left\Vert J_{2}+I_{1}^{\prime
\prime }\right\Vert _{X}+\left\Vert J_{2}+I_{1}^{\prime \prime \prime
}\right\Vert _{X}%
\end{array}%
\right) +\left\Vert I_{2}\right\Vert _{X}.
\end{equation*}%
By Lemma 1%
\begin{eqnarray*}
&&\left\Vert J_{1}^{\prime }\right\Vert _{X}\leq \frac{1}{\pi }\int_{0}^{%
\frac{\pi }{r\left( n+1\right) }}\left\Vert \psi _{\cdot }\left( t\right)
\right\Vert _{X}\left\vert \sum_{k=0}^{\infty }a_{n,k}\widetilde{D}%
_{k,1}\left( t\right) \right\vert dt \\
&\leq &\frac{1}{2\pi }\int_{0}^{\frac{\pi }{r\left( n+1\right) }}\left\Vert
\psi _{\cdot }\left( t\right) \right\Vert _{X}\sum_{k=0}^{\infty }a_{n,k}%
\frac{\pi }{t}dt\leq \frac{1}{2}\int_{0}^{\frac{\pi }{r\left( n+1\right) }}%
\frac{\omega \left( t\right) }{t}dt.
\end{eqnarray*}%
Since, by Lemma 2,%
\begin{eqnarray*}
&&\sum_{k=0}^{\infty }a_{n,k}\widetilde{D^{\circ }}_{k,1}\left( t\right)
=\sum_{k=0}^{\infty }a_{n,k}\frac{\cos \frac{\left( 2k+1\right) t}{2}}{2\sin 
\frac{t}{2}} \\
&=&\frac{1}{2\sin \frac{t}{2}}\left( \sum_{k=0}^{\infty }a_{n,k}\cos kt\cos 
\frac{t}{2}-\sum_{k=0}^{\infty }a_{n,k}\sin kt\sin \frac{t}{2}\right)
\end{eqnarray*}%
\begin{eqnarray*}
&=&\frac{\cos \frac{t}{2}}{2\sin \frac{t}{2}}\left( \sum_{k=0}^{\infty
}\left( a_{n,k}-a_{n,k+r}\right) D_{k,r}^{\circ }\left( t\right)
+\sum_{k=0}^{r-1}a_{n,k}D_{k,-r}^{\circ }\left( t\right) \right) \\
&&-\frac{1}{2}\left( -\sum_{k=0}^{\infty }\left( a_{n,k}-a_{n,k+r}\right) 
\widetilde{D^{\circ }}_{k,r}\left( t\right) -\sum_{k=0}^{r-1}a_{n,k}%
\widetilde{D^{\circ }}_{k,-r}\left( t\right) \right) ,
\end{eqnarray*}%
whence%
\begin{equation*}
\left\vert \sum_{k=0}^{\infty }a_{n,k}\widetilde{D^{\circ }}_{k,1}\left(
t\right) \right\vert \leq \frac{1}{2\left\vert \sin \frac{t}{2}\sin \frac{rt%
}{2}\right\vert }\left( A_{n,r}+\sum_{k=0}^{r-1}a_{n,k}\right) \leq \frac{1}{%
\left\vert \sin \frac{t}{2}\sin \frac{rt}{2}\right\vert }A_{n,r}.
\end{equation*}%
Hence and by Lemma 1, 
\begin{eqnarray*}
&&\left( 
\begin{array}{c}
\left\Vert J_{1}+J_{2}+I_{1}^{\prime \prime }\right\Vert _{X} \\ 
\left\Vert J_{2}+I_{1}^{\prime \prime }\right\Vert _{X}%
\end{array}%
\right) \\
&\leq &\left( 
\begin{array}{c}
\frac{1}{\pi }\sum\limits_{m=0}^{\left[ r/2\right] }\left( \int_{\frac{2m\pi 
}{r}}^{\frac{2m\pi }{r}+\frac{\pi }{r\left( n+1\right) }}+\int_{\frac{2m\pi 
}{r}+\frac{\pi }{r\left( n+1\right) }}^{\frac{2m\pi }{r}+\frac{\pi }{r}%
}\right) \left\Vert \psi _{\cdot }\left( t\right) \right\Vert _{X}\left\vert
\sum\limits_{k=0}^{\infty }a_{n,k}\widetilde{D^{\circ }}_{k,1}\left(
t\right) \right\vert dt \\ 
\frac{1}{\pi }\left( \sum\limits_{m=1}^{\left[ r/2\right] }\int_{\frac{2m\pi 
}{r}}^{\frac{2m\pi }{r}+\frac{\pi }{r\left( n+1\right) }}+\sum\limits_{m=0}^{%
\left[ r/2\right] }\int_{\frac{2m\pi }{r}+\frac{\pi }{r\left( n+1\right) }}^{%
\frac{2m\pi }{r}+\frac{\pi }{r}}\right) \left\Vert \psi _{\cdot }\left(
t\right) \right\Vert _{X}\left\vert \sum\limits_{k=0}^{\infty }a_{n,k}%
\widetilde{D^{\circ }}_{k,1}\left( t\right) \right\vert dt%
\end{array}%
\right)
\end{eqnarray*}%
and therefore%
\begin{eqnarray*}
&&\left. 
\begin{array}{c}
\left\Vert J_{1}+J_{2}+I_{1}^{\prime \prime }\right\Vert _{X} \\ 
\left\Vert J_{1}^{\prime }\right\Vert _{X}+\left\Vert J_{2}+I_{1}^{\prime
\prime }\right\Vert _{X}%
\end{array}%
\right\} \\
&\leq &\frac{1}{2}\sum_{m=0}^{\left[ r/2\right] }\int_{\frac{2m\pi }{r}}^{%
\frac{2m\pi }{r}+\frac{\pi }{r\left( n+1\right) }}\frac{O\left( \omega
\left( t\right) \right) }{t}dt+\frac{1}{\pi }\sum_{m=0}^{\left[ r/2\right]
}\int_{\frac{2m\pi }{r}+\frac{\pi }{r\left( n+1\right) }}^{\frac{2m\pi }{r}+%
\frac{\pi }{r}}\frac{O\left( \omega \left( t\right) \right) }{\left\vert
\sin \frac{t}{2}\sin \frac{rt}{2}\right\vert }A_{n,r}dt.
\end{eqnarray*}%
Using Lemmas 3, 4, with $c=b=r$, and the estimates $\left\vert \sin \frac{t}{%
2}\right\vert \geq \frac{\left\vert t\right\vert }{\pi }$ for $t\in \left[
0,\pi \right] ,$ $\left\vert \sin \frac{rt}{2}\right\vert \geq \frac{rt}{\pi 
}-2m$ for $t\in \left[ \frac{2m\pi }{r},\frac{2m\pi }{r}+\frac{\pi }{r}%
\right] ,$ where $m\in \left\{ 0,...,\left[ \ r/2\right] \right\} ,$ we
obtain 
\begin{equation*}
\left. 
\begin{array}{c}
\left\Vert J_{1}+J_{2}+I_{1}^{\prime \prime }\right\Vert _{X} \\ 
\left\Vert J_{1}^{\prime }\right\Vert _{X}+\left\Vert J_{2}+I_{1}^{\prime
\prime }\right\Vert _{X}%
\end{array}%
\right\} \leq O\left( 1\right) \left( \left[ r/2\right] +1\right) \frac{\pi 
}{r\left( n+1\right) }H\left( \frac{\pi }{n+1}\right)
\end{equation*}%
\begin{equation*}
+\sum_{k=0}^{\infty }\left\vert a_{n,k}-a_{n,k+r}\right\vert \sum_{m=0}^{ 
\left[ r/2\right] }\int_{\frac{2m\pi }{r}+\frac{\pi }{r\left( n+1\right) }}^{%
\frac{2m\pi }{r}+\frac{\pi }{r}}\frac{\omega \left( t\right) }{t\left( \frac{%
rt}{\pi }-2m\right) }dt
\end{equation*}%
\begin{equation*}
=O\left( 1\right) \frac{\pi }{n+1}H\left( \frac{\pi }{n+1}\right)
+A_{n,r}\sum_{m=0}^{\left[ r/2\right] }\int_{\frac{2m\pi }{r}+\frac{\pi }{%
r\left( n+1\right) }}^{\frac{2m\pi }{r}+\frac{\pi }{r}}\frac{O\left( \omega
\left( t\right) \right) }{\frac{rt}{\pi }\left( t-\frac{2m\pi }{r}\right) }dt
\end{equation*}%
\begin{equation*}
=O\left( 1\right) \frac{\pi }{n+1}H\left( \frac{\pi }{n+1}\right)
+A_{n,r}\sum_{m=0}^{\left[ r/2\right] }\int_{\frac{\pi }{r\left( n+1\right) }%
}^{\frac{\pi }{r}}\frac{O\left( \omega \left( t+\frac{2m\pi }{r}\right)
\right) }{\frac{rt}{\pi }\left( t+\frac{2m\pi }{r}\right) }dt
\end{equation*}%
\begin{equation*}
\leq O\left( 1\right) \left[ \frac{\pi }{n+1}H\left( \frac{\pi }{n+1}\right)
+\left( \left[ r/2\right] +1\right) \frac{2\pi }{r}A_{n,r}\int_{\frac{\pi }{%
r\left( n+1\right) }}^{\frac{\pi }{r}}\frac{\omega \left( t\right) }{t^{2}}dt%
\right]
\end{equation*}%
\begin{equation*}
=O\left( 1\right) \left[ \frac{\pi }{n+1}H\left( \frac{\pi }{n+1}\right)
+A_{n,r}H\left( \frac{\pi }{n+1}\right) \right] .
\end{equation*}%
Analogously%
\begin{equation*}
\left. 
\begin{array}{c}
\left\Vert J_{1}+J_{2}+I_{1}^{\prime \prime \prime }\right\Vert _{X} \\ 
\left\Vert J_{1}^{\prime }\right\Vert _{X}+\left\Vert J_{2}+I_{1}^{\prime
\prime \prime }\right\Vert _{X}%
\end{array}%
\right\} =O\left( 1\right) \left[ \frac{\pi }{r\left( n+1\right) }H\left( 
\frac{\pi }{n+1}\right) +\frac{2\pi }{r}A_{n,r}H\left( \frac{\pi }{n+1}%
\right) \right] .
\end{equation*}%
Similarly, by Lemma 1, Lemmas 3, 4, with $c=b=r$ and the estimates $%
\left\vert \sin \frac{t}{2}\right\vert \geq \frac{\left\vert t\right\vert }{%
\pi }$ for $t\in \left[ 0,\pi \right] ,$ $\left\vert \sin \frac{rt}{2}%
\right\vert \geq 2\left( m+1\right) -\frac{rt}{\pi }$ for $t\in \left[ \frac{%
2\left( m+1\right) \pi }{r}-\frac{\pi }{r},\frac{2\left( m+1\right) \pi }{r}-%
\frac{\pi }{r\left( n+1\right) }\right] ,$ where $m\in \left\{ 0,...,\left[
r/2\right] -1\right\} ,$ we get%
\begin{eqnarray*}
&&\left\Vert I_{2}\right\Vert _{X}\leq \frac{1}{\pi }\sum_{m=0}^{\left[ r/2%
\right] -1}\int_{\frac{2m\pi }{r}+\frac{\pi }{r}}^{\frac{2\left( m+1\right)
\pi }{r}}\left\Vert \psi _{\cdot }\left( t\right) \right\Vert _{X}\left\vert
\sum_{k=0}^{\infty }a_{n,k}\widetilde{D^{\circ }}_{k,1}\left( t\right)
\right\vert dt \\
&=&\frac{1}{\pi }\sum_{m=0}^{\left[ r/2\right] -1}\left( \int_{\frac{2\left(
m+1\right) \pi }{r}-\frac{\pi }{r}}^{\frac{2\left( m+1\right) \pi }{r}-\frac{%
\pi }{r\left( n+1\right) }}+\int_{\frac{2\left( m+1\right) \pi }{r}-\frac{%
\pi }{r\left( n+1\right) }}^{\frac{2\left( m+1\right) \pi }{r}}\right)
\left\Vert \psi _{\cdot }\left( t\right) \right\Vert _{X}\left\vert
\sum_{k=0}^{\infty }a_{n,k}\widetilde{D^{\circ }}_{k,1}\left( t\right)
\right\vert dt
\end{eqnarray*}%
\begin{equation*}
\leq \frac{1}{2}\sum_{m=0}^{\left[ r/2\right] -1}\int\limits_{\frac{2\left(
m+1\right) \pi }{r}-\frac{\pi }{r\left( n+1\right) }}^{\frac{2\left(
m+1\right) \pi }{r}}\frac{O\left( \omega \left( t\right) \right) }{t}dt+%
\frac{1}{\pi }A_{n,r}\sum_{m=0}^{\left[ r/2\right] -1}\int\limits_{\frac{%
2\left( m+1\right) \pi }{r}-\frac{\pi }{r}}^{\frac{2\left( m+1\right) \pi }{r%
}-\frac{\pi }{r\left( n+1\right) }}\frac{O\left( \omega \left( t\right)
\right) }{\left\vert \sin \frac{t}{2}\sin \frac{rt}{2}\right\vert }dt
\end{equation*}%
\begin{equation*}
\leq \frac{1}{2}\sum_{m=0}^{\left[ r/2\right] -1}\int\limits_{\frac{2\left(
m+1\right) \pi }{r}-\frac{\pi }{r\left( n+1\right) }}^{\frac{2\left(
m+1\right) \pi }{r}}\frac{O\left( \omega \left( t\right) \right) }{t}%
dt+A_{n,r}\sum_{m=0}^{\left[ r/2\right] -1}\int\limits_{\frac{2\left(
m+1\right) \pi }{r}-\frac{\pi }{r}}^{\frac{2\left( m+1\right) \pi }{r}-\frac{%
\pi }{r\left( n+1\right) }}\frac{O\left( \omega \left( t\right) \right) }{%
\frac{rt}{\pi }\left[ \frac{2\left( m+1\right) \pi }{r}-t\right] }dt
\end{equation*}%
\begin{equation*}
=\frac{1}{2}\sum_{m=0}^{\left[ r/2\right] -1}\int\limits_{\frac{2\left(
m+1\right) \pi }{r}-\frac{\pi }{r\left( n+1\right) }}^{\frac{2\left(
m+1\right) \pi }{r}}\frac{O\left( \omega \left( t\right) \right) }{t}%
dt+A_{n,r}\sum_{m=0}^{\left[ r/2\right] -1}\int\limits_{\frac{\pi }{r\left(
n+1\right) }}^{\frac{\pi }{r}}\frac{O\left( \omega \left( -t+\frac{2\left(
m+1\right) \pi }{r}\right) \right) }{\frac{r}{\pi }t\left( -t+\frac{2\left(
m+1\right) \pi }{r}\right) }dt
\end{equation*}%
\begin{equation*}
\leq \frac{1}{2}\sum_{m=0}^{\left[ r/2\right] -1}\int_{\frac{2\left(
m+1\right) \pi }{r}-\frac{\pi }{r\left( n+1\right) }}^{\frac{2\left(
m+1\right) \pi }{r}}\frac{O\left( \omega \left( t\right) \right) }{t}%
dt+A_{n,r}\left[ r/2\right] \frac{2\pi }{r}\int_{\frac{\pi }{r\left(
n+1\right) }}^{\frac{\pi }{r}}\frac{O\left( \omega \left( t\right) \right) }{%
t^{2}}dt.
\end{equation*}%
Thus%
\begin{equation*}
\left\Vert I_{2}\right\Vert _{X}=O\left( 1\right) \left[ \frac{\pi }{n+1}%
H\left( \frac{\pi }{n+1}\right) +A_{n,r}H\left( \frac{\pi }{n+1}\right) %
\right] .
\end{equation*}%
Collecting these estimates we obtain the first result.

Applying condition $\left( \ref{113}\right) $ we have 
\begin{equation*}
\left[ \left( n+1\right) \sum_{k=0}^{\infty }\left\vert
a_{n,k}-a_{n,k+r}\right\vert \right] ^{-1}=\left[ \sum_{l=0}^{n}\sum_{k=0}^{%
\infty }\left\vert a_{n,k}-a_{n,k+r}\right\vert \right] ^{-1}\leq \left[
\sum_{l=0}^{n}\sum_{k=l}^{\infty }\left\vert a_{n,k}-a_{n,k+r}\right\vert %
\right] ^{-1}
\end{equation*}%
\begin{equation*}
\leq \left[ \sum_{l=0}^{n}\left\vert \sum_{k=l}^{\infty }\left(
a_{n,k}-a_{n,k+r}\right) \right\vert \right] ^{-1}\leq \left[
\sum_{l=0}^{n}\sum_{k=l}^{r+l-1}a_{n,k}\right] ^{-1}=O\left( 1\right)
\end{equation*}%
and the second result also follows. $\blacksquare $

\subsection{Proof of Theorem 2}

Analogously, as in the proof of Theorem 1, we consider an odd $r$ and an
even $r.$ Then,%
\begin{equation*}
\left( 
\begin{array}{c}
\widetilde{T}_{n,A}^{\text{ }}f\left( x\right) -\widetilde{f}\left( x\right)
\\ 
\widetilde{T}_{n,A}^{\text{ }}f\left( x\right) -\widetilde{f}\left( x,\frac{1%
}{r}A_{n,r}\right)%
\end{array}%
\right)
\end{equation*}%
\begin{eqnarray*}
&=&\left( 
\begin{array}{c}
+ \\ 
-%
\end{array}%
\right) \frac{1}{\pi }\int\limits_{0}^{\frac{1}{r}A_{n,r}}\psi _{x}\left(
t\right) \sum_{k=0}^{\infty }a_{n,k}\left( 
\begin{array}{c}
\widetilde{D^{\circ }}_{k,1}\left( t\right) \\ 
\widetilde{D}_{k,1}\left( t\right)%
\end{array}%
\right) dt+\frac{1}{\pi }\int\limits_{\frac{1}{r}A_{n,r}}^{\frac{\pi }{r}%
}\psi _{x}\left( t\right) \sum_{k=0}^{\infty }a_{n,k}\widetilde{D^{\circ }}%
_{k,1}\left( t\right) dt \\
&&+\frac{1}{\pi }\sum_{m=1}^{\left[ r/2\right] }\int\limits_{\frac{2m\pi }{r}%
}^{\frac{2m\pi }{r}+\frac{\pi }{r}}\psi _{x}\left( t\right)
\sum_{k=0}^{\infty }a_{n,k}\widetilde{D^{\circ }}_{k,1}\left( t\right) dt+%
\frac{1}{\pi }\sum_{m=0}^{\left[ r/2\right] -1}\int\limits_{\frac{2m\pi }{r}+%
\frac{\pi }{r}}^{\frac{2\left( m+1\right) \pi }{r}}\psi _{x}\left( t\right)
\sum_{k=0}^{\infty }a_{n,k}\widetilde{D^{\circ }}_{k,1}\left( t\right) dt
\end{eqnarray*}%
or%
\begin{equation*}
\left( 
\begin{array}{c}
\widetilde{T}_{n,A}^{\text{ }}f\left( x\right) -\widetilde{f}\left( x\right)
\\ 
\widetilde{T}_{n,A}^{\text{ }}f\left( x\right) -\widetilde{f}\left( x,\frac{1%
}{r}A_{n,r}\right)%
\end{array}%
\right)
\end{equation*}%
\begin{equation*}
=\left( 
\begin{array}{c}
+ \\ 
-%
\end{array}%
\right) \frac{1}{\pi }\int\limits_{0}^{\frac{1}{r}A_{n,r}}\psi _{x}\left(
t\right) \sum_{k=0}^{\infty }a_{n,k}\left( 
\begin{array}{c}
\widetilde{D^{\circ }}_{k,1}\left( t\right) \\ 
\widetilde{D}_{k,1}\left( t\right)%
\end{array}%
\right) dt+\frac{1}{\pi }\int\limits_{\frac{1}{r}A_{n,r}}^{\frac{\pi }{r}%
}\psi _{x}\left( t\right) \sum_{k=0}^{\infty }a_{n,k}\widetilde{D^{\circ }}%
_{k,1}\left( t\right) dt
\end{equation*}%
\begin{equation*}
+\frac{1}{\pi }\sum_{m=1}^{\left[ r/2\right] -1}\int\limits_{\frac{2m\pi }{r}%
}^{\frac{2m\pi }{r}+\frac{\pi }{r}}\psi _{x}\left( t\right)
\sum_{k=0}^{\infty }a_{n,k}\widetilde{D^{\circ }}_{k,1}\left( t\right) dt+%
\frac{1}{\pi }\sum_{m=0}^{\left[ r/2\right] -1}\int\limits_{\frac{2m\pi }{r}+%
\frac{\pi }{r}}^{\frac{2\left( m+1\right) \pi }{r}}\psi _{x}\left( t\right)
\sum_{k=0}^{\infty }a_{n,k}\widetilde{D^{\circ }}_{k,1}\left( t\right) dt,
\end{equation*}%
respectively. Since $A_{n,r}\leq 2,$ therefore we can estimate our terms
analogously as in the proof of Theorem 1 with $A_{n,r}$ instead of $\frac{%
\pi }{n+1}$ and thus we obtain the desired estimate. $\blacksquare $

\subsection{Proof of Theorem 3}

Similarly, as in the proof of Theorem 1%
\begin{equation*}
\left( 
\begin{array}{c}
\left\Vert \widetilde{T}_{n,A}^{\text{ }}f\left( \cdot \right) -\widetilde{f}%
\left( \cdot \right) \right\Vert _{X} \\ 
\left\Vert \widetilde{T}_{n,A}^{\text{ }}f\left( \cdot \right) -\widetilde{f}%
\left( \cdot ,\frac{\pi }{r\left( n+1\right) }\right) \right\Vert _{X}%
\end{array}%
\right) \leq \left( 
\begin{array}{c}
\left\Vert J_{1}\right\Vert _{X} \\ 
\left\Vert J_{1}^{\prime }\right\Vert _{X}%
\end{array}%
\right) +\left\Vert J_{2}+I_{1}^{\prime \prime }\right\Vert _{X}+\left\Vert
J_{2}+I_{1}^{\prime \prime \prime }\right\Vert _{X}+\left\Vert
I_{2}\right\Vert _{X}.
\end{equation*}%
By Lemma 1 and $\left( \ref{114}\right) $ 
\begin{eqnarray*}
\left\Vert J_{1}^{\prime }\right\Vert _{X} &\leq &\frac{1}{\pi }%
\int\limits_{0}^{\frac{\pi }{r\left( n+1\right) }}\left\Vert \psi _{\cdot
}\left( t\right) \right\Vert _{X}\left\vert \sum_{k=0}^{\infty }a_{n,k}%
\widetilde{D}_{k,1}\left( t\right) \right\vert dt\leq \frac{1}{\pi }%
\sum_{k=0}^{\infty }\left( k+1\right) a_{n,k}\int\limits_{0}^{\frac{\pi }{%
r\left( n+1\right) }}\omega \left( t\right) dt \\
&=&O\left( n+1\right) \int_{0}^{\frac{\pi }{r\left( n+1\right) }}\omega
\left( t\right) dt=O\left( 1\right) \omega \left( \frac{\pi }{r\left(
n+1\right) }\right) =O\left( \omega \left( \frac{\pi }{n+1}\right) \right)
\end{eqnarray*}%
and by Lemma 1 and $\left( \ref{202}\right) $ 
\begin{eqnarray*}
\left\Vert J_{1}\right\Vert _{X} &\leq &\frac{1}{\pi }\int\limits_{0}^{\frac{%
\pi }{r\left( n+1\right) }}\left\Vert \psi _{\cdot }\left( t\right)
\right\Vert _{X}\left\vert \sum_{k=0}^{\infty }a_{n,k}\widetilde{D^{\circ }}%
_{k,1}\left( t\right) \right\vert dt\leq \frac{1}{2\pi }\int\limits_{0}^{%
\frac{\pi }{r\left( n+1\right) }}\left\Vert \psi _{\cdot }\left( t\right)
\right\Vert _{X}\sum_{k=0}^{\infty }a_{n,k}\frac{\pi }{t}dt \\
&\leq &\frac{1}{2}\int_{0}^{\frac{\pi }{r\left( n+1\right) }}\frac{\omega
\left( t\right) }{t}dt=O\left( \omega \left( \frac{\pi }{r\left( n+1\right) }%
\right) \right) =O\left( \omega \left( \frac{\pi }{n+1}\right) \right) .
\end{eqnarray*}%
Further, by the same lemmas and conditions as in the above proofs and Lemma
5, we obtain with $\kappa =\left\{ 
\begin{array}{c}
1\text{ when }r\text{ is even,} \\ 
0\text{ when }r\text{ is odd,}%
\end{array}%
\right. $ that 
\begin{eqnarray*}
&&\left\Vert J_{2}+I_{1}^{\prime \prime }\right\Vert _{X}+\left\Vert
J_{2}+I_{1}^{\prime \prime \prime }\right\Vert _{X} \\
&\leq &\frac{1}{\pi }\left( \sum_{m=1}^{\left[ r/2\right] -\kappa }\int_{%
\frac{2m\pi }{r}}^{\frac{2m\pi }{r}+\frac{\pi }{r}}+\int_{\frac{\pi }{%
r\left( n+1\right) }}^{\frac{\pi }{r}}\right) \left\Vert \psi _{\cdot
}\left( t\right) \right\Vert _{X}\left\vert \sum_{k=0}^{\infty }a_{n,k}%
\widetilde{D^{\circ }}_{k,1}\left( t\right) \right\vert dt \\
&=&\frac{1}{\pi }\left( \sum_{m=1}^{\left[ r/2\right] -\kappa }\int_{\frac{%
2m\pi }{r}}^{\frac{2m\pi }{r}+\frac{\pi }{r\left( n+1\right) }}+\sum_{m=0}^{%
\left[ r/2\right] -\kappa }\int_{\frac{2m\pi }{r}+\frac{\pi }{r\left(
n+1\right) }}^{\frac{2m\pi }{r}+\frac{\pi }{r}}\right) \left\Vert \psi
_{\cdot }\left( t\right) \right\Vert _{X}\left\vert \sum_{k=0}^{\infty
}a_{n,k}\widetilde{D^{\circ }}_{k,1}\left( t\right) \right\vert dt
\end{eqnarray*}%
\begin{equation*}
\leq \frac{1}{\pi }\sum_{m=1}^{\left[ r/2\right] -\kappa }\int\limits_{\frac{%
2m\pi }{r}}^{\frac{2m\pi }{r}+\frac{\pi }{r\left( n+1\right) }}\frac{O\left(
\omega \left( t\right) \right) }{2\left\vert \sin \frac{t}{2}\right\vert }dt+%
\frac{1}{\pi }\sum_{m=0}^{\left[ r/2\right] -\kappa }\int\limits_{\frac{%
2m\pi }{r}+\frac{\pi }{r\left( n+1\right) }}^{\frac{2m\pi }{r}+\frac{\pi }{r}%
}\frac{O\left( \omega \left( t\right) \right) }{\left\vert \sin \frac{t}{2}%
\sin \frac{rt}{2}\right\vert }A_{n,r}dt
\end{equation*}%
\begin{eqnarray*}
&\leq &\frac{1}{2}\sum_{m=1}^{\left[ r/2\right] -\kappa }\int_{\frac{2m\pi }{%
r}}^{\frac{2m\pi }{r}+\frac{\pi }{r\left( n+1\right) }}\frac{O\left( \omega
\left( t\right) \right) }{t}dt+A_{n,r}\sum_{m=0}^{\left[ r/2\right] -\kappa
}\int_{\frac{2m\pi }{r}+\frac{\pi }{r\left( n+1\right) }}^{\frac{2m\pi }{r}+%
\frac{\pi }{r}}\frac{O\left( \omega \left( t\right) \right) }{t\left( \frac{%
rt}{\pi }-2m\right) }dt \\
&\leq &\sum_{m=1}^{\left[ r/2\right] -\kappa }\frac{O\left( \omega \left( 
\frac{2m\pi }{r}\right) \right) }{\frac{2m\pi }{r}}\int_{\frac{2m\pi }{r}}^{%
\frac{2m\pi }{r}+\frac{\pi }{r\left( n+1\right) }}dt+A_{n,r}\sum_{m=0}^{%
\left[ r/2\right] -\kappa }\int_{\frac{\pi }{r\left( n+1\right) }}^{\frac{%
\pi }{r}}\frac{O\left( \omega \left( t+\frac{2m\pi }{r}\right) \right) }{%
\frac{rt}{\pi }\left( t+\frac{2m\pi }{r}\right) }dt
\end{eqnarray*}%
\begin{eqnarray*}
&\leq &2\sum_{m=1}^{\left[ r/2\right] -\kappa }\frac{O\left( \omega \left( 
\frac{2\pi }{r}\right) \right) }{\frac{2\pi }{r}}\frac{\pi }{r\left(
n+1\right) }+\frac{2\pi }{r}\left( \left[ r/2\right] +1\right) A_{n,r}\int_{%
\frac{\pi }{r\left( n+1\right) }}^{\frac{\pi }{r}}\frac{O\left( \omega
\left( t\right) \right) }{t^{2}}dt \\
&=&O\left( 1\right) \left[ \omega \left( \frac{\pi }{n+1}\right)
+A_{n,r}H\left( \frac{\pi }{n+1}\right) \right]
\end{eqnarray*}%
and%
\begin{equation*}
\left\Vert I_{2}\right\Vert _{X}\leq \frac{1}{\pi }\sum_{m=0}^{\left[ r/2%
\right] -1}\left( \int_{\frac{2\left( m+1\right) \pi }{r}-\frac{\pi }{r}}^{%
\frac{2\left( m+1\right) \pi }{r}-\frac{\pi }{r\left( n+1\right) }}+\int_{%
\frac{2\left( m+1\right) \pi }{r}-\frac{\pi }{r\left( n+1\right) }}^{\frac{%
2\left( m+1\right) \pi }{r}}\right) \left\Vert \psi _{\cdot }\left( t\right)
\right\Vert _{X}\left\vert \sum_{k=0}^{\infty }a_{n,k}\widetilde{D^{\circ }}%
_{k,1}\left( t\right) \right\vert dt
\end{equation*}%
\begin{equation*}
\leq \frac{1}{\pi }\sum_{m=0}^{\left[ r/2\right] -1}\int\limits_{\frac{%
2\left( m+1\right) \pi }{r}-\frac{\pi }{r}}^{\frac{2\left( m+1\right) \pi }{r%
}-\frac{\pi }{r\left( n+1\right) }}\frac{O\left( \omega \left( t\right)
\right) }{\left\vert \sin \frac{t}{2}\sin \frac{rt}{2}\right\vert }A_{n,r}dt+%
\frac{1}{2}\sum_{m=0}^{\left[ r/2\right] -1}\int\limits_{\frac{2\left(
m+1\right) \pi }{r}-\frac{\pi }{r\left( n+1\right) }}^{\frac{2\left(
m+1\right) \pi }{r}}\frac{O\left( \omega \left( t\right) \right) }{t}dt
\end{equation*}%
\begin{equation*}
\leq \sum_{m=0}^{\left[ r/2\right] -1}A_{n,r}\int\limits_{\frac{2\left(
m+1\right) \pi }{r}-\frac{\pi }{r}}^{\frac{2\left( m+1\right) \pi }{r}-\frac{%
\pi }{r\left( n+1\right) }}\frac{O\left( \omega \left( t\right) \right) }{%
\frac{rt}{\pi }\left[ \frac{2\left( m+1\right) \pi }{r}-t\right] }%
dt+\sum_{m=0}^{\left[ r/2\right] -1}\frac{O\left( \omega \left( \frac{%
2\left( m+1\right) \pi }{r}-\frac{\pi }{r\left( n+1\right) }\right) \right) 
}{r\left( n+1\right) \left( \frac{2\left( m+1\right) \pi }{r}-\frac{\pi }{%
r\left( n+1\right) }\right) }
\end{equation*}%
\begin{eqnarray*}
&\leq &\left[ r/2\right] \left[ \frac{2\pi }{r}A_{n,r}\int_{\frac{\pi }{%
r\left( n+1\right) }}^{\frac{\pi }{r}}\frac{O\left( \omega \left( t\right)
\right) }{t^{2}}dt+2\frac{O\left( \omega \left( \frac{\pi }{r}\right)
\right) }{\frac{\pi }{r}}\frac{\pi }{r\left( n+1\right) }\right] \\
&=&O\left( 1\right) \left[ A_{n,r}H\left( \frac{\pi }{n+1}\right) +\omega
\left( \frac{\pi }{n+1}\right) \right] .
\end{eqnarray*}%
Thus our proof is complete. $\blacksquare $

\subsection{Proof of Theorem 4}

Let as above%
\begin{equation*}
\left( 
\begin{array}{c}
\left\Vert \widetilde{T}_{n,A}^{\text{ }}f\left( \cdot \right) -\widetilde{f}%
\left( \cdot \right) \right\Vert _{X} \\ 
\left\Vert \widetilde{T}_{n,A}^{\text{ }}f\left( \cdot \right) -\widetilde{f}%
\left( \cdot ,\frac{\pi }{r\left( n+1\right) }\right) \right\Vert _{X}%
\end{array}%
\right) \leq \left( 
\begin{array}{c}
\left\Vert J_{1}\right\Vert _{X} \\ 
\left\Vert J_{1}^{\prime }\right\Vert _{X}%
\end{array}%
\right) +\left\Vert J_{2}+I_{1}^{\prime \prime }\right\Vert _{X}+\left\Vert
J_{2}+I_{1}^{\prime \prime \prime }\right\Vert _{X}+\left\Vert
I_{2}\right\Vert _{X}\text{ },
\end{equation*}%
\begin{eqnarray*}
&&\left\Vert J_{1}\right\Vert _{X}\leq \frac{1}{\pi }\int\limits_{0}^{\frac{%
\pi }{r\left( n+1\right) }}\left\Vert \psi _{\cdot }\left( t\right)
\right\Vert _{X}\left\vert \sum_{k=0}^{\infty }a_{n,k}\widetilde{D^{\circ }}%
_{k,1}\left( t\right) \right\vert dt\leq \frac{1}{2}\sum_{k=0}^{\infty
}a_{n,k}\int\limits_{0}^{\frac{\pi }{r\left( n+1\right) }}\frac{\widetilde{%
\omega }_{2}\left( f,t\right) _{X}}{t}dt \\
&=&O\left( 1\right) \widetilde{\omega }_{2}\left( f,\frac{\pi }{r\left(
n+1\right) }\right) _{X}=O\left( \widetilde{\omega }_{2}\left( f,\frac{\pi }{%
n+1}\right) _{X}\right)
\end{eqnarray*}%
and%
\begin{eqnarray*}
&&\left\Vert J_{1}^{\prime }\right\Vert _{X} \\
&\leq &\frac{1}{\pi }\int\limits_{0}^{\frac{\pi }{r\left( n+1\right) }%
}\left\Vert \psi _{\cdot }\left( t\right) \right\Vert _{X}\left\vert
\sum_{k=0}^{\infty }a_{n,k}\widetilde{D}_{k,1}\left( t\right) \right\vert
dt\leq \frac{1}{2\pi }\int\limits_{0}^{\frac{\pi }{r\left( n+1\right) }%
}\left\Vert \psi _{\cdot }\left( t\right) \right\Vert _{X}\sum_{k=0}^{\infty
}a_{n,k}\left( k+1\right) dt \\
&\leq &O\left( n+1\right) \int_{0}^{\frac{\pi }{r\left( n+1\right) }}%
\widetilde{\omega }_{2}\left( f,t\right) _{X}dt=O\left( \widetilde{\omega }%
_{2}\left( f,\frac{\pi }{n+1}\right) _{X}\right) ,
\end{eqnarray*}%
Further, taking $\tau _{m}^{1}=\left[ \frac{\pi }{rt-2m\pi }\right] $ and $%
\tau =\left[ \frac{\pi }{rt}\right] $, using Lemma 5, we obtain with $\kappa
=\left\{ 
\begin{array}{c}
1\text{ when }r\text{ is even,} \\ 
0\text{ when }r\text{ is odd,}%
\end{array}%
\right. $ that 
\begin{eqnarray*}
&&\left\Vert J_{2}+I_{1}^{\prime \prime }\right\Vert _{X}+\left\Vert
J_{2}+I_{1}^{\prime \prime \prime }\right\Vert _{X} \\
&\leq &\frac{1}{\pi }\left( \sum_{m=1}^{\left[ r/2\right] -\kappa }\int_{%
\frac{2m\pi }{r}}^{\frac{2m\pi }{r}+\frac{\pi }{r}}+\int_{\frac{\pi }{%
r\left( n+1\right) }}^{\frac{\pi }{r}}\right) \left\Vert \psi _{\cdot
}\left( t\right) \right\Vert _{X}\left\vert \sum_{k=0}^{\infty }a_{n,k}%
\widetilde{D^{\circ }}_{k,1}\left( t\right) \right\vert dt \\
&=&\frac{1}{\pi }\left( \sum_{m=1}^{\left[ r/2\right] -\kappa }\int_{\frac{%
2m\pi }{r}}^{\frac{2m\pi }{r}+\frac{\pi }{r\left( n+1\right) }}+\sum_{m=0}^{%
\left[ r/2\right] -\kappa }\int_{\frac{2m\pi }{r}+\frac{\pi }{r\left(
n+1\right) }}^{\frac{2m\pi }{r}+\frac{\pi }{r}}\right) \left\Vert \psi
_{\cdot }\left( t\right) \right\Vert _{X}\left\vert \sum_{k=0}^{\infty
}a_{n,k}\widetilde{D^{\circ }}_{k,1}\left( t\right) \right\vert dt
\end{eqnarray*}%
\begin{eqnarray*}
&\leq &\frac{1}{\pi }\sum_{m=1}^{\left[ r/2\right] -\kappa }\int_{\frac{%
2m\pi }{r}}^{\frac{2m\pi }{r}+\frac{\pi }{r\left( n+1\right) }}\frac{%
\widetilde{\omega }_{2}\left( f,t\right) _{X}}{2\left\vert \sin \frac{t}{2}%
\right\vert }\sum_{k=0}^{\infty }a_{n,k}dt \\
&&+\frac{1}{\pi }\sum_{m=0}^{\left[ r/2\right] -\kappa }\int\limits_{\frac{%
2m\pi }{r}+\frac{\pi }{r\left( n+1\right) }}^{\frac{2m\pi }{r}+\frac{\pi }{r}%
}\left( \frac{\widetilde{\omega }_{2}\left( f,t\right) _{X}}{2\left\vert
\sin \frac{t}{2}\right\vert }\sum_{k=0}^{\tau _{m}^{1}}a_{n,k}+\frac{%
\widetilde{\omega }_{2}\left( f,t\right) _{X}}{\left\vert \sin \frac{t}{2}%
\sin \frac{rt}{2}\right\vert }\sum_{k=\tau _{m}^{1}}^{\infty }\left\vert
a_{n,k}-a_{n,k+r}\right\vert \right) dt
\end{eqnarray*}%
\begin{eqnarray*}
&\leq &\frac{1}{2}\sum_{m=1}^{\left[ r/2\right] -\kappa }\int\limits_{\frac{%
2m\pi }{r}}^{\frac{2m\pi }{r}+\frac{\pi }{r\left( n+1\right) }}\frac{%
\widetilde{\omega }_{2}\left( f,t\right) _{X}}{t}dt+\left( \left[ r/2\right]
+1\right) \int_{\frac{\pi }{r\left( n+1\right) }}^{\frac{\pi }{r}}\frac{%
\widetilde{\omega }_{2}\left( f,t\right) _{X}}{t}\sum_{k=0}^{\tau }a_{n,k}dt
\\
&&+\frac{2\pi }{r}\left( \left[ r/2\right] +1\right) \int_{\frac{\pi }{%
r\left( n+1\right) }}^{\frac{\pi }{r}}\frac{\widetilde{\omega }_{2}\left(
f,t\right) _{X}}{t^{2}}\sum_{k=\tau }^{\infty }\left\vert
a_{n,k}-a_{n,k+r}\right\vert dt
\end{eqnarray*}%
\begin{eqnarray*}
&\leq &\sum_{m=1}^{\left[ r/2\right] -\kappa }\frac{\widetilde{\omega }%
_{2}\left( f,\frac{2m\pi }{r}\right) _{X}}{\frac{2m\pi }{r}}\int\limits_{%
\frac{2m\pi }{r}}^{\frac{2m\pi }{r}+\frac{\pi }{r\left( n+1\right) }%
}dt+\left( \left[ r/2\right] +1\right) \sum_{\mu =1}^{n}\int\limits_{\mu
}^{\mu +1}\frac{\widetilde{\omega }_{2}\left( f,\frac{\pi }{rt}\right) _{X}}{%
\frac{\pi }{rt}}\sum_{k=0}^{\left[ t\right] }a_{n,k}\frac{\pi dt}{rt^{2}} \\
&&+\frac{2\pi }{r}\left( \left[ r/2\right] +1\right) \sum_{\mu
=1}^{n}\int\limits_{\mu }^{\mu +1}\frac{\widetilde{\omega }_{2}\left( f,%
\frac{\pi }{rt}\right) _{X}}{\left( \frac{\pi }{rt}\right) ^{2}}\sum_{k=%
\left[ t\right] }^{\infty }\left\vert a_{n,k}-a_{n,k+r}\right\vert \frac{\pi
dt}{rt^{2}}
\end{eqnarray*}%
\begin{eqnarray*}
&\leq &O\left( 1\right) \widetilde{\omega }_{2}\left( f,\frac{\pi }{n+1}%
\right) _{X}+O\left( 1\right) \sum_{\mu =1}^{n}\frac{\widetilde{\omega }%
_{2}\left( f,\frac{\pi }{\mu }\right) _{X}}{\mu }\sum_{k=0}^{\mu +1}a_{n,k}
\\
&&+O\left( 1\right) \sum_{\mu =1}^{n}\widetilde{\omega }_{2}\left( f,\frac{%
\pi }{\mu }\right) _{X}\sum_{k=\mu }^{\infty }\left\vert
a_{n,k}-a_{n,k+r}\right\vert .
\end{eqnarray*}%
Next, taking $\tau _{m}^{2}=\left[ \frac{\pi }{-rt+2\left( m+1\right) \pi }%
\right] $ , we obtain 
\begin{eqnarray*}
\left\Vert I_{2}\right\Vert _{X} &\leq &\frac{1}{\pi }\sum_{m=0}^{\left[ r/2%
\right] -1}\int_{\frac{2m\pi }{r}+\frac{\pi }{r}}^{\frac{2\left( m+1\right)
\pi }{r}}\left\Vert \psi _{\cdot }\left( t\right) \right\Vert _{X}\left\vert
\sum_{k=0}^{\infty }a_{n,k}\widetilde{D^{\circ }}_{k,1}\left( t\right)
\right\vert dt \\
&\leq &\frac{1}{\pi }\sum_{m=0}^{\left[ r/2\right] -1}\left( \int\limits_{%
\frac{2\left( m+1\right) \pi }{r}-\frac{\pi }{r}}^{\frac{2\left( m+1\right)
\pi }{r}-\frac{\pi }{r\left( n+1\right) }}+\int\limits_{\frac{2\left(
m+1\right) \pi }{r}-\frac{\pi }{r\left( n+1\right) }}^{\frac{2\left(
m+1\right) \pi }{r}}\right) \left\Vert \psi _{\cdot }\left( t\right)
\right\Vert _{X}\left\vert \sum_{k=0}^{\infty }a_{n,k}\widetilde{D^{\circ }}%
_{k,1}\left( t\right) \right\vert dt
\end{eqnarray*}%
\begin{eqnarray*}
&\leq &\frac{1}{\pi }\sum_{m=0}^{\left[ r/2\right] -1}\int\limits_{\frac{%
2\left( m+1\right) \pi }{r}-\frac{\pi }{r}}^{\frac{2\left( m+1\right) \pi }{r%
}-\frac{\pi }{r\left( n+1\right) }}\left( \frac{\widetilde{\omega }%
_{2}\left( f,t\right) _{X}}{2\left\vert \sin \frac{t}{2}\right\vert }%
\sum_{k=0}^{\tau _{m}^{2}}a_{n,k}\right. +\left. \frac{\widetilde{\omega }%
_{2}\left( f,t\right) _{X}}{\left\vert \sin \frac{t}{2}\sin \frac{rt}{2}%
\right\vert }\sum_{k=\tau _{m}^{2}}^{\infty }\left\vert
a_{n,k}-a_{n,k+r}\right\vert dt\right) \\
&&+\frac{1}{\pi }\sum_{m=0}^{\left[ r/2\right] -1}\int_{\frac{2\left(
m+1\right) \pi }{r}-\frac{\pi }{r\left( n+1\right) }}^{\frac{2\left(
m+1\right) \pi }{r}}\frac{\widetilde{\omega }_{2}\left( f,t\right) _{X}}{%
2\left\vert \sin \frac{t}{2}\right\vert }\sum_{k=0}^{\infty }a_{n,k}dt
\end{eqnarray*}%
\begin{equation*}
\leq \frac{1}{2}\sum_{m=0}^{\left[ r/2\right] -1}\int_{\frac{\pi }{r\left(
n+1\right) }}^{\frac{\pi }{r}}\frac{\widetilde{\omega }_{2}\left( f,-t+\frac{%
2\left( m+1\right) \pi }{r}\right) _{X}}{-t+\frac{2\left( m+1\right) \pi }{r}%
}\sum_{k=0}^{\tau }a_{n,k}dt
\end{equation*}%
\begin{equation*}
+\sum_{m=0}^{\left[ r/2\right] -1}\int\limits_{\frac{\pi }{r\left(
n+1\right) }}^{\frac{\pi }{r}}\frac{\widetilde{\omega }_{2}\left( f,-t+\frac{%
2\left( m+1\right) \pi }{r}\right) _{X}}{\frac{r}{\pi }t\left( -t+\frac{%
2\left( m+1\right) \pi }{r}\right) }\sum_{k=\tau }^{\infty }\left\vert
a_{n,k}-a_{n,k+r}\right\vert dt+\frac{1}{2}\sum_{m=0}^{\left[ r/2\right]
-1}\int\limits_{\frac{2\left( m+1\right) \pi }{r}-\frac{\pi }{r\left(
n+1\right) }}^{\frac{2\left( m+1\right) \pi }{r}}\frac{\widetilde{\omega }%
_{2}\left( f,t\right) _{X}}{t}dt
\end{equation*}%
\begin{equation*}
\leq \left[ r/2\right] \int\limits_{\frac{\pi }{r\left( n+1\right) }}^{\frac{%
\pi }{r}}\frac{\widetilde{\omega }_{2}\left( f,t\right) _{X}}{t}%
\sum_{k=0}^{\tau }a_{n,k}+\frac{2\pi }{r}\left[ r/2\right] \int\limits_{%
\frac{\pi }{r\left( n+1\right) }}^{\frac{\pi }{r}}\frac{\widetilde{\omega }%
_{2}\left( f,t\right) _{X}}{t^{2}}\sum_{k=\tau }^{\infty }\left\vert
a_{n,k}-a_{n,k+r}\right\vert dt
\end{equation*}%
\begin{equation*}
+\sum_{m=0}^{\left[ r/2\right] -1}\frac{\widetilde{\omega }_{2}\left( f,%
\frac{2\left( m+1\right) \pi }{r}-\frac{\pi }{r\left( n+1\right) }\right)
_{X}}{\frac{2\left( m+1\right) \pi }{r}-\frac{\pi }{r\left( n+1\right) }}%
\int_{\frac{2\left( m+1\right) \pi }{r}-\frac{\pi }{r\left( n+1\right) }}^{%
\frac{2\left( m+1\right) \pi }{r}}dt
\end{equation*}%
\begin{eqnarray*}
&\leq &O\left( 1\right) \sum_{\mu =1}^{n}\frac{\widetilde{\omega }_{2}\left(
f,\frac{\pi }{\mu }\right) _{X}}{\mu }\sum_{k=0}^{\mu +1}a_{n,k}+O\left(
1\right) \sum_{\mu =1}^{n}\widetilde{\omega }_{2}\left( f,\frac{\pi }{\mu }%
\right) _{X}\sum_{k=\mu }^{\infty }\left\vert a_{n,k}-a_{n,k+r}\right\vert \\
&&+O\left( 1\right) \widetilde{\omega }_{2}\left( f,\frac{\pi }{n+1}\right)
_{X}.
\end{eqnarray*}

Thus the result follows. $\blacksquare $

\subsection{Proof of Corollary 1}

Theorem 4 implies that%
\begin{eqnarray*}
&&\left. 
\begin{array}{c}
\left\Vert \widetilde{T}_{n,A}^{\text{ }}f\left( \cdot \right) -\widetilde{f}%
\left( \cdot \right) \right\Vert _{X} \\ 
\left\Vert \widetilde{T}_{n,A}^{\text{ }}f\left( \cdot \right) -\widetilde{f}%
\left( \cdot ,\frac{\pi }{r\left( n+1\right) }\right) \right\Vert _{X}%
\end{array}%
\right\} =O\left( \widetilde{\omega }_{2}\left( f,\frac{\pi }{n+1}\right)
_{X}+\sum_{\mu =1}^{n}\frac{\widetilde{\omega }_{2}\left( f,\frac{\pi }{\mu }%
\right) _{X}}{\mu }\sum_{k=0}^{\mu +1}a_{n,k}\right. \\
&&+\left. \sum_{\mu =1}^{n}\widetilde{\omega }_{2}\left( f,\frac{\pi }{\mu }%
\right) _{X}\sum_{k=\mu }^{\infty }\left\vert a_{n,k}-a_{n,k+r}\right\vert
\right) .
\end{eqnarray*}%
Since $\left( \ref{200}\right) $ 
\begin{equation*}
\sum_{\mu =1}^{n}\widetilde{\omega }_{2}\left( f,\frac{\pi }{\mu }\right)
_{X}\sum_{k=\mu }^{\infty }\left\vert a_{n,k}-a_{n,k+r}\right\vert =O\left(
1\right) \sum_{\mu =1}^{n}\widetilde{\omega }_{2}\left( f,\frac{\pi }{\mu }%
\right) _{X}\left( \sum_{k=\mu /c}^{\mu -1}+\sum_{k=\mu
}^{n}+\sum_{k=n+1}^{\infty }\right) \frac{a_{n,k}}{k}
\end{equation*}%
\begin{equation*}
\leq O\left( 1\right) \sum_{\mu =1}^{n}\widetilde{\omega }_{2}\left( f,\frac{%
\pi }{\mu }\right) _{X}\left( \sum_{k=\mu /c}^{\mu -1}\frac{a_{n,k}}{k}%
\right) +O\left( 1\right) \sum_{\mu =1}^{n}\widetilde{\omega }_{2}\left( f,%
\frac{\pi }{\mu }\right) _{X}\left( \sum_{k=\mu }^{n}+\sum_{k=n+1}^{\infty
}\right) \frac{a_{n,k}}{k}
\end{equation*}%
\begin{equation*}
\leq O\left( 1\right) c\sum_{\mu =1}^{n}\frac{\widetilde{\omega }_{2}\left(
f,\frac{\pi }{\mu }\right) _{X}}{\mu }\sum_{k=0}^{\mu +1}a_{n,k}+O\left(
1\right) \sum_{\mu =1}^{n}\widetilde{\omega }_{2}\left( f,\frac{\pi }{\mu }%
\right) _{X}\left( \sum_{k=\mu }^{n}+\sum_{k=n+1}^{\infty }\right) \frac{%
a_{n,k}}{k}
\end{equation*}%
one has%
\begin{eqnarray*}
&&\left. 
\begin{array}{c}
\left\Vert \widetilde{T}_{n,A}^{\text{ }}f\left( \cdot \right) -\widetilde{f}%
\left( \cdot \right) \right\Vert _{X} \\ 
\left\Vert \widetilde{T}_{n,A}^{\text{ }}f\left( \cdot \right) -\widetilde{f}%
\left( \cdot ,\frac{\pi }{r\left( n+1\right) }\right) \right\Vert _{X}%
\end{array}%
\right\} \\
&=&O\left( 1\right) \widetilde{\omega }_{2}\left( f,\frac{\pi }{n+1}\right)
_{X}+O\left( 1\right) \left( 1+c\right) \sum_{\mu =1}^{n}\frac{\widetilde{%
\omega }_{2}\left( f,\frac{\pi }{\mu }\right) _{X}}{\mu }\sum_{k=0}^{\mu
+1}a_{n,k} \\
&&+O\left( 1\right) \sum_{\mu =1}^{n}\widetilde{\omega }_{2}\left( f,\frac{%
\pi }{\mu }\right) _{X}\sum_{k=\mu }^{n}\frac{a_{n,k}}{k}+O\left( 1\right)
\sum_{k=n+1}^{\infty }\frac{a_{n,k}}{k}\sum_{\mu =1}^{n}\widetilde{\omega }%
_{2}\left( f,\frac{\pi }{\mu }\right) _{X}
\end{eqnarray*}%
\begin{eqnarray*}
&\leq &O\left( 1\right) \widetilde{\omega }_{2}\left( f,\frac{\pi }{n+1}%
\right) _{X}+O\left( 1\right) \left( 1+c\right) \left\{ \sum_{\mu =1}^{n}%
\frac{\widetilde{\omega }_{2}\left( f,\frac{\pi }{\mu }\right) _{X}}{\mu }%
\sum_{k=0}^{\mu -1}a_{n,k}\right. \\
&&+\left. 2\sum_{\mu =1}^{n}\frac{\widetilde{\omega }_{2}\left( f,\frac{\pi 
}{\mu }\right) _{X}}{\mu +1}a_{n,\mu +1}+\sum_{\mu =1}^{n}\frac{\widetilde{%
\omega }_{2}\left( f,\frac{\pi }{\mu }\right) _{X}}{\mu }a_{n,\mu }\right\}
\\
&&+O\left( 1\right) \sum_{\mu =1}^{n}\widetilde{\omega }_{2}\left( f,\frac{%
\pi }{\mu }\right) _{X}\sum_{k=\mu }^{n}\frac{a_{n,k}}{k}+O\left( 1\right)
\sum_{k=n+1}^{\infty }\frac{a_{n,k}}{k}\sum_{\mu =1}^{n}\widetilde{\omega }%
_{2}\left( f,\frac{\pi }{\mu }\right) _{X}
\end{eqnarray*}%
\begin{eqnarray*}
&\leq &O\left( 1\right) \widetilde{\omega }_{2}\left( f,\frac{\pi }{n+1}%
\right) _{X}+O\left( 1\right) \left( 1+c\right) \sum_{\mu =0}^{n}\frac{%
\widetilde{\omega }_{2}\left( f,\frac{\pi }{\mu +1}\right) _{X}}{\mu +1}%
\sum_{k=0}^{\mu }a_{n,k} \\
&&+O\left( 1\right) \left[ 3\left( 1+c\right) +1\right] \sum_{\mu =1}^{n}%
\widetilde{\omega }_{2}\left( f,\frac{\pi }{\mu }\right) _{X}\sum_{k=\mu
}^{n}\frac{a_{n,k}}{k} \\
&&+O\left( 1\right) \sum_{k=n+1}^{\infty }\frac{a_{n,k}}{k}\sum_{\mu =1}^{n}%
\widetilde{\omega }_{2}\left( f,\frac{\pi }{\mu }\right) _{X}
\end{eqnarray*}%
\begin{eqnarray*}
&=&O\left( 1\right) \widetilde{\omega }_{2}\left( f,\frac{\pi }{n+1}\right)
_{X}+O\left( 1\right) \sum_{k=0}^{n}a_{n,k}\sum_{\mu =k}^{n}\frac{\widetilde{%
\omega }_{2}\left( f,\frac{\pi }{\mu +1}\right) _{X}}{\mu +1} \\
&&+O\left( 1\right) \sum_{k=1}^{n}\frac{a_{n,k}}{k}\sum_{\mu =1}^{k}%
\widetilde{\omega }_{2}\left( f,\frac{\pi }{\mu }\right) _{X}+O\left(
1\right) \sum_{k=n+1}^{\infty }\frac{a_{n,k}}{k}\sum_{\mu =1}^{n}\widetilde{%
\omega }_{2}\left( f,\frac{\pi }{\mu }\right) _{X}.
\end{eqnarray*}%
If $\left( \ref{111}\right) $ and $\left( \ref{112}\right) $ hold, then 
\begin{equation*}
\widetilde{\omega }_{2}\left( f,\frac{\pi }{n+1}\right) _{X}\leq \frac{1}{n}%
\sum_{\mu =1}^{n}\widetilde{\omega }_{2}\left( f,\frac{\pi }{\mu }\right)
_{X}\leq O\left( 1\right) 4\pi \frac{H\left( \frac{\pi }{n+1}\right) }{n+1},
\end{equation*}%
\begin{equation*}
\sum_{\mu =k}^{n}\frac{\widetilde{\omega }_{2}\left( f,\frac{\pi }{\mu +1}%
\right) _{X}}{\mu +1}\leq 8\int_{\frac{\pi }{n+2}}^{\frac{\pi }{k+1}}\frac{%
\widetilde{\omega }_{2}\left( f,t\right) _{X}}{t}dt\leq O\left( 1\right)
8\pi \frac{H\left( \frac{\pi }{k+1}\right) }{k+1}
\end{equation*}%
and therefore%
\begin{eqnarray*}
&&\left. 
\begin{array}{c}
\left\Vert \widetilde{T}_{n,A}^{\text{ }}f\left( \cdot \right) -\widetilde{f}%
\left( \cdot \right) \right\Vert _{X} \\ 
\left\Vert \widetilde{T}_{n,A}^{\text{ }}f\left( \cdot \right) -\widetilde{f}%
\left( \cdot ,\frac{\pi }{r\left( n+1\right) }\right) \right\Vert _{X}%
\end{array}%
\right\} \\
&=&O\left( \frac{H\left( \frac{\pi }{n+1}\right) }{n+1}\right) +O\left(
\sum_{k=1}^{n}a_{n,k}\frac{H\left( \frac{\pi }{k+1}\right) }{k+1}\right)
+O\left( H\left( \frac{\pi }{n+1}\right) \sum_{k=n+1}^{\infty }\frac{a_{n,k}%
}{k}\right) .
\end{eqnarray*}%
Since 
\begin{equation*}
\sum_{k=n+1}^{\infty }\frac{a_{n,k}}{k}\leq \frac{1}{n+1}\sum_{k=n+1}^{%
\infty }a_{n,k}\leq \frac{1}{n+1}
\end{equation*}%
the result follows. $\blacksquare $

\end{document}